\documentclass[11pt]{amsart}
\usepackage{color}
\usepackage{epsfig}
\usepackage{psfrag, graphics}
\usepackage{amsthm}
\usepackage{amssymb}
\usepackage{amsmath}
\usepackage{amscd}
\usepackage{graphicx}
\newcommand{\eqdef}{\;{:=}\;}

\newtheorem {Theorem}   {Theorem}

\numberwithin{Theorem}{section}

\newtheorem {Lemma}[Theorem]    {Lemma}
\newtheorem {Proposition}[Theorem]{Proposition}
\theoremstyle{definition}

\theoremstyle{remark}
\newtheorem{Remark}[Theorem]{Remark}
\newtheorem{Example}[Theorem]{Example}

\expandafter\chardef\csname pre amssym.def
at\endcsname=\the\catcode`\@ \catcode`\@=11
\def\undefine#1{\let#1\undefined}
\def\newsymbol#1#2#3#4#5{\let\next@\relax
 \ifnum#2=\@ne\let\next@\msafam@\else
 \ifnum#2=\tw@\let\next@\msbfam@\fi\fi
 \mathchardef#1="#3\next@#4#5}
\def\mathhexbox@#1#2#3{\relax
 \ifmmode\mathpalette{}{\m@th\mathchar"#1#2#3}%
 \else\leavevmode\hbox{$\m@th\mathchar"#1#2#3$}\fi}
\def\hexnumber@#1{\ifcase#1 0\or 1\or 2\or 3\or 4\or 5\or 6\or 7\or 8\or
 9\or A\or B\or C\or D\or E\or F\fi}

\font\teneufm=eufm10 \font\seveneufm=eufm7 \font\fiveeufm=eufm5
\newfam\eufmfam
\textfont\eufmfam=\teneufm \scriptfont\eufmfam=\seveneufm
\scriptscriptfont\eufmfam=\fiveeufm

\catcode`\@=\csname pre amssym.def at\endcsname

\newcommand{\AC}{{\mathcal A}}

\newcommand{\CC}{{\mathcal C}}

\newcommand{\FF}{{\mathcal F}}

\newcommand{\HH}{{\mathcal H}}

\newcommand{\JJ}{{\mathcal J}}

\newcommand{\LL}{{\mathcal L}}

\newcommand{\NN}{{\mathcal N}}

\newcommand{\PP}{{\mathcal P}}

\newcommand{\RR}{{\mathcal R}}

\newcommand{\UU}{{\mathcal U}}
\newcommand{\VV}{{\mathcal V}}

\def    \C      {{\mathbb C}}
\def    \R      {{\mathbb R}}
\def    \Z      {{\mathbb Z}}

\newcommand{\id}{{\mathit id}}

\def    \om       {\omega}
\def    \eps      {\epsilon}
\def    \ssminus  {\smallsetminus}
\def    \p        {\partial}
\def    \codim    {\mathit{codim}\,}
\def    \dim    {\mathit{dim}\,}
\def    \length   {\operatorname{length}}

\def    \ev       {\operatorname{ev}}

\def    \Crit      {\operatorname{Crit}}

\def    \Hh     {\mathbf{H}}

\def    \12      {\frac{1}{2}}
\def    \Vol      {\operatorname{Vol}}
\def    \sr        { \Phi_R} 
\def    \sl         { \Phi_L}
\def    \srl         { \Phi_{L,k}} 
\def    \sll        { \Phi_{R,k}} 
\def    \ind      {\operatorname{ind}}

\def    \H   {\operatorname{H}}  
\def    \CF   {\operatorname{CF}} 
\def    \CM   {\operatorname{CM}} 
\def    \CZ   {\operatorname{\mu_{\scriptscriptstyle{CZ}}}} 

\begin{document}






\title[Displacement energy of coisotropic submanifolds]{Displacement energy of coisotropic submanifolds and Hofer's geometry}

\author[Ely Kerman]{Ely Kerman}
\address{Department of Mathematics, University of Illinois at Urbana-Champaign, Urbana, IL 61801, USA }
\email{ekerman@math.uiuc.edu}

\date{\today}

\thanks{This research was partially supported by NSF Grant DMS-0405994 and a grant from the Campus Research
Board of the University of Illinois at Urbana-Champaign.}


\bigskip

\begin{abstract} We prove that the displacement energy of a stable 
coisotropic submanifold is bounded away from zero if the ambient symplectic 
manifold is closed, rational and satisfies a mild topological condition. 
\end{abstract}

\maketitle

\section{Introduction and Results}

There is positive lower bound for the amount of energy it takes to displace  a closed Lagrangian submanifold of a tame symplectic manifold.  
In particular, every time-dependent function on a symplectic manifold determines a unique Hamiltonian diffeomorphism, and if this 
diffeomorphism displaces a closed Lagrangian submanifold, then the Hofer norm of the function is bounded away from zero by a constant which depends only on the Lagrangian submanifold and the 
ambient symplectic manifold.
This fundamental fact in symplectic topology was first established for rational Lagrangian submanifolds by Polterovich in \cite{po}, and was later extended to general Lagrangians by Chekanov in \cite{ch}. 
Among other things, it implies the nondegeneracy of the Hofer metric on the group of Hamiltonian diffeomorphisms of a tame symplectic manifold.

Recently, Ginzburg proved that there is also a positive lower bound for the amount of energy required to displace certain coisotropic submanifolds. More precisely, in \cite{gi:coisotropic} it is shown that 
the displacement energy of a stable coisotropic submanifold of a tame, wide  and symplectically aspherical symplectic manifold is bounded away from zero. In the present paper, we extend this coisotropic intersection phenomenon to symplectic manifolds which admit symplectic spheres. The proof utilizes the Floer theoretic methods developed in \cite{ke2}, as well as the applications of these methods to the study of Hamiltonian paths which are length minimizing  with respect to the Hofer metric. 

There is currently no version of Floer theory for the intersection theory of a general coisotropic submanifold and its image under a Hamiltonian diffeomorphism. However, one can study 
the symplectic properties of a coisotropic submanifold  using the Hamiltonian Floer homology of  functions which are supported in (normal) neighborhoods of it. This indirect approach, which goes back to the pioneering work of Viterbo from \cite{vi:torus}, requires a few compromises. 

The first compromise involves the submanifolds. To get useful normal neighborhoods, we restrict our attention to \emph{stable} coisotropic submanifolds.
This notion  was introduced by Bolle in \cite{bo1,bo2}\footnote{Bolle refers to such manifolds as being of \emph{almost contact type.}} and is defined as follows. 
Let $(M, \om)$ be a symplectic manifold of dimension $2m$ and let $N$ be a closed coisotropic submanifold of $M$ with codimension $k$. Then $N$ is said to be {\bf stable} 
if there are one-forms $\alpha_1, \dots, \alpha_k$ on $N$ such that 
the form
$$
\alpha_1 \wedge \dots \wedge \alpha_k \wedge (\om|_{\scriptscriptstyle{N}})^{m-k}
$$
does not vanish on $N$, and $\ker d\alpha_j \supset \ker \om|_{\scriptscriptstyle{N}}$ for $j = 1, \dots , k$. Examples of stable coisotropic submanifolds include Lagrangian tori and contact hypersurfaces. 
The stability condition is also closed under products.  For more details, the reader is referred  to \cite{bo1,bo2, gi:coisotropic}.

In using Hamiltonian Floer homology to study the symplectic topology of  a coisotropic submanifold, one  also needs to recognize nontrivial $1$-periodic orbits using only the symplectic action and/or the Conley-Zehnder index. 
This requires further compromise concerning the ambient symplectic manifolds, $(M,\om)$,  we consider. In  \cite{gi:coisotropic}, the symplectic manifolds are assumed to be symplectically aspherical. That is, for every 
class $A \in \pi_2(M)$  it is assumed that  $\om(A) =0 = c_1(A)$,  where the notations $\om(A)$ and $c_1(A)$ refer to the evaluations of the cohomology classes on the elements of $\H_2(M;\R)$ and $\H_2(M;\Z)$ determined by $A$.
 With this assumption, the action and index of a periodic orbit are single-valued and any $1$-periodic orbit with 
sufficiently large action (greater than $\|H\|^+$ as defined below) must be 
nonconstant. 

Here, we allow for the existence of nontrivial symplectic spheres and so the action
and index may be multi-valued. To distinguish nonconstant periodic orbits we will assume that the quantity
$$
r(M,\om) = \inf_{A \in \pi_2(M)} \left\{ |\om(A)| \mid |\om(A)|>0\right\}.
$$
is positive.\footnote{We use the convention that the infimum over the empty set is equal to $\infty$.}  
A symplectic manifold with $r(M, \om)>0$ is said to be {\bf rational}. We will also assume that $(M, \om)$ satisfies 
the topological assumption
\begin{equation}\label{link}
\om(A)=0  \Longrightarrow c_1(A) \geq 0 \text{  for all $A$ in $\pi_2(M)$}.
\end{equation}
Finally, we restrict ourselves, in this work, to the case when $M$ is closed. We expect that the 
methods developed here are also applicable to symplectic manifolds which are open or have 
convex boundaries.

Before stating the main result, we first recall the definition of the displacement energy. 
Let  $C^{\infty}(S^1 \times M)$ be the space of smooth time-periodic functions on $M$, 
where $S^1=\R/\Z$ is the circle parameterized by $t \in [0,1]$.  The Hofer norm of a function  $H$ in $C^{\infty}(S^1 \times M)$ is
defined as
\begin{eqnarray*}
 \|H\|&=& \int_0^1 \max_{p \in M}H_t(p) \,\,dt-  \int_0^1 \min_{p \in M} H_t(p) \,\,dt,
\end{eqnarray*}
where $H_t(\cdot) =  H(t, \cdot)$.
One can also associate to $H$ its Hamiltonian vector
$X_H$ via the equation
$$
\om(X_H, \cdot) = -dH_t(\cdot).
$$  
The time-$t$ flow of this vector field, also referred to as the Hamiltonian flow of $H$, is denoted by  $\phi^t_H$ and is defined for all $t \in [0,1]$.
The group of Hamiltonian diffeomorphisms consists of all the time-1 maps, $\phi^1_H$, obtained in this way.

The {\bf displacement energy} of a subset $U$ of $M$ is defined as
$$
e(U) = \inf_{H \in C^{\infty}(S^1 \times M)}\{\|H\| \mid \phi^1_H(U) \cap U=\emptyset\},
$$
the minimum variation of a function which generates a Hamiltonian diffeomorphism that
 moves $U$ off of itself.

Our main result is the following:
 
\begin{Theorem}\label{coisotropic}
Let $N$ be a stable coisotropic submanifold of a closed and rational symplectic 
manifold satisfying \eqref{link}. There is a positive constant $\Delta>0$ such 
that $e(N) \geq \Delta$.
\end{Theorem}

Of course, one starts with the assumption that $N$ can be displaced by some
Hamiltonian diffeomorphism, i.e., $e(N)<\infty$.
This has deep implications for the Hamiltonian flows supported near $N$.
In turn, these flows can be used to probe the geometry of $N$. It is this interaction 
between the displacability of $N$ and its geometry, which leads to the proof of   
Theorem \ref{coisotropic}.


The primary difference between the proof of Theorem \ref{coisotropic} and the proof of the main result in \cite{gi:coisotropic} is the contribution coming from  Floer theory. 
In \cite{gi:coisotropic}, both the action filtration and action selector are used 
to prove the existence of a Floer trajectory whose energy yields the crucial estimate for the displacement energy, (Proposition 5.1 of \cite{gi:coisotropic}). 
For a  rational symplectic manifold, the action filtration  and selector can not be used in the same manner. Instead 
we use the Floer theoretic techniques which were developed in \cite{ke2} to study the length 
minimizing properties of Hamiltonian paths.  These tools allow us to detect a perturbed holomorphic cylinder in Proposition \ref{keycap} whose 
energy recovers the crucial estimate.

\begin{Remark}
 Another approach to studying coisotropic intersections is to consider general leaf-wise intersections under Hamiltonian diffeomorphisms, \cite{mo, eh, ho1}. The most recent work in this direction 
is \cite{dr}, where Dragnev establishes the existence of leaf-wise intersections
for a stable coisotropic submanifold $N$ of $\R^{2n}$, and its image under any Hamiltonian diffeomorphism with energy less than the Floer-Hofer capacity of $N$.
\end{Remark}

\subsection{Organization}
The proof of Theorem \ref{coisotropic} is described in the next section, 
assuming the contribution from Floer theory, Proposition \ref{keycap}. 
In the third section, we recall the required Floer theory methods and 
applications from \cite{ke2}. The 
proof of Proposition \ref{keycap} is then contained in the final section
of the paper.

\subsection{Acknowledgments.}
The author would like to thank Peter Albers and Viktor Ginzburg for their helpful 
comments.

\section{The proof of Theorem \ref{coisotropic} (modulo Proposition \ref{keycap}).}

Before presenting the proof of Theorem \ref{coisotropic} in \S \ref{outline}, 
we discuss some preliminary notions and results. 

\subsection{Properties of stable coisotropic submanifolds}
We begin by recalling some useful implications of the stability assumption. The proofs of these results
can be found in \cite{bo1,bo2, gi:coisotropic}. 

Let $N$ be a stable coisotropic submanifold of codimension $k$ in a symplectic manifold $(M, \om)$
of dimension $2m$.
We then have the following normal neighborhood result.
\begin{Proposition}[\cite{bo1,bo2}]
For sufficiently small $r>0$ there is a neighborhood of $N$ in $(M, \om)$ which is 
symplectomorphic to 
$$
U_r = \{ (q,p) \in N \times \R^k \mid |p| < r\}
$$
equipped with the symplectic form 
$$
 \om|_N+ \sum_{j=1}^k d(p_j \pi^* \alpha_j).
$$
\end{Proposition}
\noindent Here, $|p|$ denotes the standard norm of $p=(p_1,\dots,p_k) \in \R^k$, and 
$\pi \colon U_r \to N$ is the obvious projection.

Recall that  the characteristic foliation $\FF$ of $N$ is determined by the 
integrable distribution $\ker \om|_{\scriptscriptstyle{N}}$.
The normal form above implies that  for each manifold $N_p = N \times p$  with $|p|<r$
we have $\om|_{\scriptscriptstyle{N_p}} =\om|_N$. Hence, each of the  $N_p$
in the tubular neighborhood $U_r$ is a coisotropic submanifold with the same characteristic foliation.

The relevant Hamiltonian dynamical system is the following leaf-wise geodesic flow on the tubular neighborhood $U_r$ of $N$. 
\begin{Proposition}[\cite{bo1,bo2,gi:coisotropic}]
\label{metric}
The Hamiltonian flow of the function $\12 |p|^2$ on the normal neighborhood $U_r$ is the 
leaf-wise geodesic flow of the leaf-wise metric $\sum_{j=1}^k (\alpha_j)^2$ on $\FF$. Moreover, 
this metric is leaf-wise flat.
\end{Proposition}

This implies that a nonconstant periodic orbit $x$ of the flow of $\12 |p|^2$ 
corresponds to a closed geodesic $\gamma$  contained on a leaf of $\FF$ in $N$.  
The fact that the leaf-wise metric is flat implies that this geodesic is noncontractible 
within its leaf.

For any  closed curve $\gamma$ contained in a leaf of $\FF$, set
$$
\delta (\gamma) = \sum_{j=1}^k \left| \int_{\gamma} \alpha_j \right|.
$$

\begin{Lemma}[\cite{bo1,bo2,gi:coisotropic}]
\label{delta}
There is a constant $\delta_N>0$ such that $$\delta(\gamma) \geq \delta_N $$ for every nontrivial closed 
geodesic $\gamma$ of the leaf-wise  metric $\sum_{j=1}^k (\alpha_j)^2.$
\end{Lemma}

\subsection{Hofer's length functional}
A function $H$ in $C^{\infty}(S^1\times M)$ is said to be normalized if 
$$
\int_M H_t \, \om^m = 0
$$
for all $t$ in $[0,1]$. The space of normalized functions is denoted by $C^{\infty}_0(S^1 \times M)$.
For every path of Hamiltonian diffeomorphisms, $\psi_t$, there is a unique $H$ in 
 $C^{\infty}_0(S^1 \times M)$ such that $\phi^t_H \circ \psi_0 = \psi_t.$ 
 Following \cite{ho1}, this time-dependent generating function
is used to define the Hofer length of the path $\psi_t$ by 
\begin{eqnarray*}
\length (\psi_t) &=& \|H\| \\ 
{}&=& \int_0^1 \max_M H_t \,\,dt-  \int_0^1 \min_M H_t \,\,dt\\
{} &=& \|H\|^+ + \|H\|^-
\end{eqnarray*}
The quantities $\|H\|^+$ and $\|H\|^-$ provide different measures of $\psi_t$ called the positive and negative Hofer lengths, respectively. The positive Hofer length will play an important role in the proof of Theorem \ref{coisotropic}.


\subsection{Right asymptotic spanning discs}\label{asd}
A  spanning disc for a loop $y \colon S^1 \to M$ is a smooth map $w$ from the
unit disc in $\C$ to $M$ such that $w(e^{2\pi i t}) =y(t)$.
A {\bf  right asymptotic spanning disc} for the 
loop $y$ is a smooth map $v \colon \R \times S^1 \to M$ such that 
\begin{itemize}
\item there is a sequence $s^-_j  \to -\infty$ for which
$$\lim_{j \to \infty} v(s^-_j,t) = y(t);$$
\item there is a  sequence $s_j^+ \to  +\infty$ for which $v(s_j^+,t)$ converges 
to a constant map $t \mapsto p$  for some point $p \in M.$
\end{itemize}
Here, convergence is with respect to the smooth topology on $C^{\infty}(S^1,M)$.

We will detect right asymptotic spanning discs for $1$-periodic orbits of the Hamiltonian flow
of a function $H$ in $C^{\infty}(S^1\times M)$. They will  be constructed using a smooth $(\R \times S^1)$-family of $\om$-compatible almost complex structures, $J_s$,  which is independent of $s \in \R$ for $|s|$ sufficiently large. This last  auxiliary structure is used to define the {\bf energy} 
of $v$ by 
$$ E(v) =  \int_0^1 \int_{-\infty}^{+\infty} \om \big( \p_s v, J_s(v)\p_s v\big)\,ds\,dt.$$
For each integer $j$, we will also consider the quantities
$$
E^j(v)=\int_0^1 \int_{s_j^-}^{-s_j^-} \om \big( \p_s v, J_s(v)\p_sv \big)\,ds\,dt.
$$
and 
$$
\AC^j_{H}(v) = \int_0^1 H(t,v(s_j^-,t))\, dt + \int_0^1 \int_{s_j^-}^{-s_j^-} \om \big( \p_s v, \p_tv \big)\,ds\,dt. 
$$

\subsection{The proof of Theorem \ref{coisotropic}}\label{outline} 

We may assume that  
\begin{equation}
\label{ass1}
3e(N)< r(M, \om),
\end{equation}
otherwise we are done. We will prove the following result which clearly implies Theorem \ref{coisotropic}.
\begin{Theorem}\label{delta}
 For sufficiently small values of  $r>0$, there is a positive constant $\Delta>0$, independent of $r$, such that $e(U_{\scriptscriptstyle{ r}})> \Delta$.
\end{Theorem}
By \eqref{ass1}, for all sufficiently small values of $r>0$ we have 
\begin{equation}\label{assumption r}
3e(U_r)<r(M,\om)
\end{equation}
Fix an $R>0$ for which this inequality holds. Henceforth, we will consider only 
neighborhoods $U_r$ for $0<r<R/2$.

In order to relate the assumption that $N$ can be displaced by a Hamiltonian diffeomorphism to 
the properties of the flow from Proposition \ref{metric}, we 
reparameterize this flow so that it extends to a global flow on $M$ which is supported in $U_r$. 

Let $\nu  \colon [0,+\infty) \to [0, +\infty)$ be a smooth function with the following properties
\begin{itemize}
  \item $\nu (0) = A$ on $[0, r/3]$;
  \item $\nu '< 0$ on $(r/3, 2r/3)$;
  \item $\nu =-B$ on $[2r/3 ,+\infty)$.
\end{itemize}
Here, $A$ and $B$ are positive constants. We then  define the function 
$$H_r(q,p)=
\begin{cases}
 \nu(|p|) & \text{when $(q,p)$ is in $U_r$}, \\
  -B    & \text{otherwise}.
\end{cases}
$$

The Hamiltonian flow of $H_r$ is trivial way from $U_r$, and inside of $U_r$ it is a reparameterization of the geodesic flow from Proposition \ref{metric}.
Clearly, $\|H_r\|^+ = A$ and $\|H_r\|^- =B$.
We choose the constant $A$ so that  
\begin{equation*}\label{A}
2e(U_r) < A < 3e(U_r).
\end{equation*}
We then choose a constant $B$ satisfying 
\begin{equation*}
0 < B <  A \frac{\Vol(U_r)}{\Vol(M \ssminus U_r)},
\end{equation*}
so that $H_r$ is normalized. Further restricting $R$, if necessary, we may also assume that 
\begin{equation}
\label{Hr bounds}
2e(U_r) < A+B =\|H_r\| < 3 e(U_r).
\end{equation}

The following technical result is proved in the final section of the paper using the methods developed in \cite{ke2}.
The existence of the map $v$ described below is implied by the fact that the Hamiltonian path generated 
by $H_r$ does not minimize the positive Hofer length in its homotopy class (see \S \ref{curve shortening}).
 
\begin{Proposition}\label{keycap}
For the function $H_r$ above, there is an $\epsilon >0$, a family of almost complex structures $J_s$ as in \S \ref{asd},  and a {\bf nonconstant} $1$-periodic orbit $y$ of $H_r$ 
with a right asymptotic spanning disc $v$ such that 
\begin{equation}\label{key1}
-B +E^j(v) \leq \AC^j_{H_r}(v) \leq A -\epsilon
\end{equation}
for all $j$. 
Moreover, $v$ is a solution of 
\begin{equation}\label{key2}
\p_sv + J_s(v)(\p_t v - X_{\widetilde{H_s}}(v))=0,
\end{equation}
where $\widetilde{H_s}$ is either the function $H_r$ or the function $(\eta(-s)-1)B + \eta (-s)H_r$
for a
smooth nondecreasing function $\eta(s)$ which equals zero for $s \leq -1$ and equals one for $s \geq 1.$
\end{Proposition}

The following inequality for the energy of the map $v$ detected in Proposition \ref{keycap}, 
is easily derived from the work of Bolle and Ginzburg. We include a 
proof for the sake of completeness.

\begin{Lemma}[\cite{bo1,bo2,gi:coisotropic}]\label{constant2}
There is a constant $c_R>0$ such that for the periodic orbit $y$ and the asymptotic right spanning disc $v$
of Proposition \ref{keycap}, we have 
$$E(v) >c_R \cdot \delta(\pi(y)).$$
\end{Lemma}

\begin{proof}
Let $\hat{f} \colon [0,R) \to \R$ be a smooth nonincreasing function which is equal to one on $[0,R/2)$ and is equal to zero near
$R$. Let $f$ be the function which equals $\hat{f}(|p|)$ in $U_R$ and vanishes outside of $U_R$.

For the one-forms $\sigma_i = f\pi^* \alpha_i$, we have 
\begin{equation}
\label{vanish}
i_{X_{H_r}} d \sigma_i  = 0.
\end{equation}
In particular, away from $U_r$ we have $X_{H_r}=0$.  Within $U_r$, $\sigma_i=\pi^*\alpha_i$ and  $X_{H_r}$ is a reparameterization of the 
leaf-wise geodesic flow from Proposition \ref{metric}. Hence, 
$$
i_{X_{H_r}} d \sigma_i =  i_{X_{H_r}} \pi^* d \alpha_i =0
$$
since the forms $d \alpha_i$ vanish on the leaves of $\FF$.
It follows from \eqref{vanish} that $i_{X_{\widetilde{H_s}}} d \sigma_i = 0$ for both of the possible  functions $\widetilde{H_s}$ from Proposition \ref{keycap}.

The $s$-norm of a tangent vector $X \in T_p M$ is defined to be 
$\|X\|_s = \om(X,J_s X)$. Since $J_s$ does not depend on 
$s$ when $|s|$ is large, we can find constants $c_i>0$ such that 
$$|d\sigma_i(X,Y)|\leq c_i \|X\|_s \cdot\|Y\|_s$$
for any pair of tangent vectors $X$, $Y \in T_pM$ and every $s \in \R$. 

For the asymptotic right spanning disc $v$
detected in Proposition \ref{keycap}, we then have 
\begin{eqnarray*}
  E(v) &=&  \int_0^1 \int_{-\infty}^{+\infty} \om \left( \p_s v, J_s(v)\p_sv \right)\,ds\,dt \\
 {} &=& \int_0^1 \int_{-\infty}^{+\infty} \| \p_s v \|_s \cdot \|\p_t v  - X_{\widetilde{H_s}}\|_s \, ds\,dt\\
 {} &\geq & c_i^{-1} \int_0^1 \int_{-\infty}^{+\infty} \left|d \sigma_i\left(\p_s v,\p_t v - X_{\widetilde{H_s}}\right) \right| \,ds\, dt\\
  {} &\geq & c_i^{-1} \lim_{j \to \infty} \int_0^1  \int_{s_j^-}^{s_j^+}\left|d \sigma_i\left(\p_s v,\p_t v \right) \right|\,ds\, dt \\
 {} &\geq& c_i^{-1} \lim_{j \to \infty} \left| \int_0^1 \int_{s_j^-}^{s_j^+}d \sigma_i \left( \p_s v,\p_t v \right)\,ds\, dt  \right| \\
 {} &=& c_i^{-1} \lim_{ j \to \infty}\left| \int_{v(s^+_j,\cdot)}  \sigma_i -  \int_{v(s_j^-,\cdot)}  \sigma_i \right|\\
  {} &=& c_i^{-1}\left| \int_{y}  \sigma_i \right|.
\end{eqnarray*}
Since $y$ is nonconstant, it is contained in $U_r$ where $\sigma_i = \pi^* \alpha_i$. The inequality 
above then implies that 
$$
E(v) \geq c_i^{-1}\left| \int_{y}  \sigma_i \right| = c_i^{-1}\left| \int_{\pi(y)}  \alpha_i \right|.
$$
Setting $c_R = \frac{1}{k}\min\{c^{-1}_i\}$, we are done.

\end{proof}

We can now complete the proof of Theorem \ref{coisotropic}.
By Proposition \ref{keycap}, we have 
$$
-B+E^j(v) \leq \AC^j_{H_r}(v)\leq A-\epsilon.
$$
for all $j$. Taking the limit as $j \to \infty$,  yields
\begin{equation}
\label{A+B}
\|H_r\| = A+B > E(v).
\end{equation}
Together with inequality \eqref{Hr bounds} and Lemmas \ref{constant2} and \ref{delta}, this implies
that
$$
e(U_r) > \|H_r\|/3 >  E(v)/3 > c_R \cdot \delta(\pi(y)) /3 \geq c_R \cdot \delta_N/3.
$$
Setting $\Delta = c_R \cdot \delta_N /3$, the proof of Theorem \ref{delta}, 
and hence Theorem \ref{coisotropic}, will be  complete once we prove Proposition \ref{keycap}.








\section{Floer caps and chain isomorphisms in Morse homology}

Throughout this section $H$ will be a normalized function in $C_0^{\infty}(S^1 \times M)$
whose contractible periodic orbits with period equal to one are nondegenerate.
This (finite) set of $1$-periodic orbits will be denoted by $\PP(H)$. 

Let $\JJ(M,\om)$ be the space of smooth  almost complex structures on $M$
which are compatible with $\om$, and let $\JJ_{S^1}(M,\om)$ denote
the space of smooth $S^1$-families of elements in $\JJ(M,\om)$.
Fixing a $J$ in  $\JJ_{S^1}(M,\om)$, we refer to $(H,J)$ as our {\bf Hamiltonian data}.
  

\subsection{Homotopy triples and Floer caps}
A smooth $\R$-family, $F_s$, of functions in $C^{\infty}(S^1 \times M)$ or elements of $\JJ_{S^1}(M,\om)$ 
is called a {\bf compact homotopy}
from $F^-$ to $F^+$, if there is a constant $\lambda >0$ such that 
$F_s = F^-$ for $s \leq -\lambda $, and $F_s = F^+$ for $s\geq \lambda $.
Any such $\lambda$ will be referred to as a horizon of the compact homotopy.

A {\bf homotopy triple} for the pair $(H,J)$ is a collection $$\HH = (H_s,K_s, J_s),$$ where $H_s$ is a 
compact homotopy from a constant function $c$ to $H$, 
$K_s$ is a compact homotopy from the zero function to itself, and $J_s$ is a compact homotopy  in $\JJ_{S^1}(M,\om)$ 
from  some $J^-$  to $J$. Although the constant $c$ is an important part of the homotopy triple, it will be suppressed from the notation for simplicity.

For a homotopy triple $\HH=(H_s,K_s,J_s)$, we consider smooth maps $u$ from the 
infinite cylinder $\R \times S^1$ to $M$ which satisfy the following
equation
\begin{equation}\label{left-section}
    \partial_s u- X_{K_s}(u)+ J_s(u)(\partial_tu - X_{H_s}(u))=0.
\end{equation}
The energy of a solution $u$ of \eqref{left-section} is defined as
$$
  E(u) = \int_0^1 \int_{-\infty}^{+\infty}\om \left( \p_su -X_{K_s}(u), J_s(u) (\p_su -X_{K_s}(u))\right) \,ds \, dt.
$$
If the energy of $u$ is finite, then it follows from standard arguments that
$$ u(+\infty) \eqdef \lim_{s \to \infty} u(s,t) = x(t)$$
for some $1$-periodic orbit
 $x \in\PP(H)$. 
The assumption of finite energy also implies that 
$$
u(-\infty) \eqdef \lim_{s \to -\infty} u(s,t) = p
$$
for some point $p \in M$.

The set of {\bf left Floer caps} of $x \in \PP(H)$ with respect to $\HH$ is 
$$\LL(x;\HH)= \left\{\ u \in C^{\infty}(\R \times S^1,M) \mid u
\text{ satisfies \eqref{left-section} },\,E(u)< \infty,\, u(\infty)=x \right\}.
$$
It is clear from the asymptotic behavior described above, that each left Floer cap $u$ in $\LL(x;\HH)$
determines  a unique homotopy class of spanning discs for $x$ and
hence a well-defined Conley-Zehnder index $\CZ(x,u)$.
This index is normalized here so that
if $x(t)=p$ is a constant $1$-periodic orbit of a $C^2$-small Morse function and $u(z)=p$ is the constant spanning disc,
then $\CZ(x,u) = \ind(p) -m$, where $\ind(p)$ is the Morse index of $p$.
The action of $x$ with
respect to $u$ is defined by
$$
\AC_H(x,u) = \int_0^1 H(t,x(t))\,dt - \int_0^1 \int_{-\infty}^{+\infty} \om(\p_s u, \p_t u)\,ds\,dt.
$$


Given any map of the form $F(s,\cdot)$ for $s \in \R$, we set $$\overleftarrow{F}(s, \cdot)= F(-s, \cdot).$$
For a homotopy triple $\HH=(H_s,K_s,J_s)$ we will also consider maps $v \colon
\R \times S^1 \to M$ which satisfy the equation
\begin{equation}\label{right-section}
    \p_sv + X_{\overleftarrow{K_s}}(v)+ \overleftarrow{J_s}(v)(\p_tv - X_{\overleftarrow{H_s}}(v))=0.
\end{equation}
The energy of such a map is defined by 
\begin{equation*}
E(v) =  \int_0^1 \int_{-\infty}^{+\infty} \om \Big( \p_sv +X_{\overleftarrow{K_s}}(v), \overleftarrow{J_s}(v)(\p_sv +X_{\overleftarrow{K_s}}(v) \Big) \, ds \, dt .
\end{equation*}
If a map $v$ satisfies \eqref{right-section} and has finite energy, 
then $v(+\infty)$ is a point in $M$ and 
$v(-\infty)$ is a $1$-periodic orbit of $H$. 
In particular, such a $v$ is a right asymptotic spanning
disc, as defined in \S \ref{asd}.

The space of {\bf right
Floer caps} for $x \in \PP(H)$ is defined by
$$\RR(x;\HH)=\left\{ v  \in C^{\infty}(\R \times S^1,M) \mid v \text{ satisfies \eqref{right-section}} ,\,E(v)< \infty,\,
v(-\infty)=x \right\}.$$
For each $v$ in $\RR(x;\HH)$, one can consider the map 
$$ \overleftarrow{v}(s,t)=v(-s,t)$$ and 
define the Conley-Zehnder index $\CZ(x,\overleftarrow{v})$ and the action
$$
\AC_H(x,\overleftarrow{v}) = \int_0^1 H(t,x(t))\,dt - \int_0^1 \int_{-\infty}^{+\infty} \om(\p_s \overleftarrow{v}, \p_t \overleftarrow{v})\,ds\,dt.
$$

\subsection{Curvature}
The  curvature of a homotopy triple $\HH=(H_s,K_s,J_s)$ is the function on $\R \times S^1 \times M$
defined by
$$\kappa(\HH) = \p_sH_s - \p_tK_s + \{H_s,K_s\}.$$
Here we use the convention $\{H,K\}= \om(X_K,X_H\}.$ 
The curvature relates the 
energy and action of solutions of equations \eqref{left-section} and \eqref{right-section}
as follows.

Given a solution $u$ of  \eqref{left-section}, set
\begin{equation*}
\label{}
E_a^b(u) = \int_0^1\int_a^b \om \Big(\p_s u -X_{K_s}(u), J_s(u)( \p_s u - X_{K_s}(u))\Big) \, ds \, dt 
\end{equation*}
and 
\begin{equation*}
\label{ }
\AC_a^b(u) = \int_0^1H_b(t,u(b,t))\,dt - \int_0^1 \int_a^b \om(\p_su, \p_tu) \,ds\,dt.
\end{equation*}
We then derive the following identity from \eqref{left-section} 
\begin{equation}
\label{action energy left}
E_a^b(u)= \int_0^1H_a(t,u(a,t)))\,dt  -  \AC_a^b(u) + \int_0^1 \int_a^b \kappa(\HH)(s,t,v(s,t))\,ds\,dt.
\end{equation}

For a solution $v$ of \eqref{right-section},
the corresponding map $\overleftarrow{v}(s,t) = v(-s,t)$ satisfies
\begin{equation}
\label{right-section-back}
-\p_s\overleftarrow{v} +  X_{K_s}(\overleftarrow{v})+ J_s(\overleftarrow{v})(\p_t\overleftarrow{v} - X_{H_s}(\overleftarrow{v}))=0.
\end{equation}
Setting
\begin{equation*}
\label{ }
\overleftarrow{E}_a^b(v) = \int_0^1\int_a^b \om \Big(\p_s\overleftarrow{v} - X_{K_s}(\overleftarrow{v}), J_s(\overleftarrow{v})( \p_s \overleftarrow{v} -X_{K_s}(\overleftarrow{v}))\Big) \, ds \, dt 
\end{equation*}
and 
\begin{equation*}
\label{reverse action}
\overleftarrow{\AC}_a^b(v) = \int_0^1H_b(t,\overleftarrow{v}(b,t))\,dt - \int_0^1 \int_a^b \om(\p_s\overleftarrow{v}, \p_t\overleftarrow{v}) \,ds\,dt
\end{equation*}
equation \eqref{right-section-back} then yields 
\begin{equation}
\label{action energy right}
\overleftarrow{E}_a^b(v) = \overleftarrow{\AC}_a^b(v) - \int_0^1H_a(t,\overleftarrow{v}(a,t))\,dt 
-\int_0^1 \int_a^b \kappa(\HH)(s,t,\overleftarrow{v}(s,t))\,ds\,dt
\end{equation}

Finally, we note that if $\lambda$ is a horizon for $H_s$ and if  $A \leq a \leq -\lambda$ and $B \geq b$, then 
\begin{equation}\label{action difference}
\overleftarrow{\AC}_A^B(v) - \overleftarrow{\AC}_a^b(v) \geq \int_0^1 \int_{[A,a] \cup [b,B]}\kappa(\HH)(s,t,\overleftarrow{v})\,ds \,dt.
\end{equation}

\begin{Remark}
\label{new terms}
In this notation, we have $$E(v) =\overleftarrow{E}^{+\infty}_{-\infty}(v)\,\,\, \text{ and } \,\,\,\AC_H(x,\overleftarrow{v}) = \overleftarrow{\AC}^{+\infty}_{-\infty}(v).$$
Moreover, for $j$ sufficiently large, we have  
$$E^j(v) = \overleftarrow{E}^{-s_j^-}_{s_j^-}(v) \,\,\, \text{ and } \,\,\, \AC^j_H(v)= \overleftarrow{\AC}^{-s_j^-}_{s_j^-}(v)$$ for the quantities  appearing in Proposition \ref{keycap}.
\end{Remark}

The positive and negative norms of the curvature are defined by
$$
|||\kappa(\HH)|||^+ = \int_{\R \times S^1} \max_{p\in
M}\kappa(\HH)\,\,ds \, dt,
$$
and
$$
|||\kappa(\HH)|||^- = -\int_{\R \times S^1}  \min_{p \in
M}\kappa(\HH)\,\,ds\,dt.
$$

\begin{Example}\label{ex:linear}
Let $\eta \colon \R \to [0,1]$ be a smooth nondecreasing function 
such that $\eta(s) =0 $ for $s \leq -1$ and $\eta(s)=1$ for $s \geq 1$.
This function will be fixed throughout this paper. 
A {\bf linear homotopy triple} for $(H,J)$ is a homotopy triple of the form 
$
\overline{\HH}=\left(\overline{H}_s, 0, J_s \right)
$
where $$\overline{H}_s=(\eta(s)-1)\|H\|^- + \eta(s)H.$$
The constant for  $\overline{\HH}$ is $c=-\|H\|^-.$ 
The curvature of $\overline{\HH}$ is 
$$
\kappa(\overline{\HH})= \dot{\eta}(s)(H +\|H\|^-),
$$
which is clearly negative. We also have 
$$
|||\kappa(\overline{\HH})|||^+ = \int_{S^1} \max_{p \in M} (\|H\|^- +H(t,p)) \, dt = \|H\| ,
$$
and
$$
|||\kappa(\overline{\HH})|||^- = -\int_{S^1} \min_{p \in M} (\|H\|^- +H(t,p)) \, dt =0.
$$ 
\end{Example}

Any function $G$ which generates 
a path that is homotopic to $\phi^t_H$, relative its endpoints,  
can be used to construct a useful homotopy triple for $H$. The following basic result in this direction is a simple consequence 
of Propositions 2.6 and 2.7 from \cite{ke2}.

\begin{Proposition}\label{good caps}
Let $H$ be function in $C^{\infty}_0(S^1 \times M)$. For any $G$ in $C^{\infty}_0(S^1 \times M)$ 
whose Hamiltonian path $\phi^t_G$ is 
homotopic to $\phi^t_H$ relative its endpoints,
there is a family of almost complex structures $J$ in $\JJ_{S^1}(M, \om)$ 
and a homotopy triple $\HH_G$ for $(H,J)$ such that 
$$
|||\kappa(\HH_G)|||^+ + c \leq \|G\|^+.
$$
Here, $c$ is the constant appearing in the homotopy triple $\HH_G$.
\end{Proposition}

\subsection{Cap data and central orbits}

For the pair $(H,J)$, we fix a pair of homotopy triples
$$
\Hh=(\HH_L,\HH_R).
$$
This will be  referred to as a choice of {\bf cap data}. The norm of the curvature of the cap data $\Hh$ is defined as
$$
|||\kappa(\Hh)||| = |||\kappa(\HH_R)|||^- + |||\kappa(\HH_L)|||^+.
$$

A periodic orbit $x \in \PP(H)$ is said to be {\bf{central}} for the cap data $\Hh$, if there is a pair $(u,v) \in \LL(x;\HH_L) \times
 \RR(x;\HH_R)$ such that $$[u \#v]=0.$$ Here, $u\#v$ denotes
the obvious concatenation of the maps, and $[u \#v]$ is the element of $\pi_2(M)$
determined by $u\#v$. We will refer to $(u,v)$ above as a
{\bf central pair} of Floer caps for $x$. 

For a period orbit $x$ in $\PP(H)$ and Floer caps  
$u \in \LL(x;\HH_L)$ and $v \in \RR(x;\HH_R)$, equations \eqref{action energy left} 
and \eqref{action energy right} yield  
\begin{equation}\label{energy-left}
    0 \leq E(u) =  c_L - \AC_H(x,u)+
\int_0^1 \int_{-\infty}^{+\infty} \kappa(\HH_L)(s,t,u(s,t))\,ds \, dt,
\end{equation}
and
\begin{equation}\label{energy-right}
    0 \leq E(v) = \AC_H(x,\overleftarrow{v}) -c_R - \int_0^1 \int_{-\infty}^{+\infty} \kappa(\HH_R)(s,t,\overleftarrow{v}(s,t))\,ds \, dt,
\end{equation}
where $c_L$ and $c_R$ are the constants for $\HH_L$ and $\HH_R$, respectively.

If $(u,v)$ is a central pair of Floer caps for $x$ with respect to $\Hh$, then 
$\AC_H(x,u)=\AC_H(x,\overleftarrow{v})$ and 
\eqref{energy-left} and \eqref{energy-right} imply that
\begin{equation}\label{action-bounds}
    -|||\kappa(\HH_R)|||^- +c_R\leq \AC_H(x,\overleftarrow{v}) = \AC_H(x,u) \leq  |||\kappa(\HH_L)|||^+ + c_L.
\end{equation}
For a central pair $(u,v)$ one also obtains from \eqref{energy-left}, \eqref{energy-right} and \eqref{action-bounds}, the uniform energy bounds
\begin{equation}\label{uniform-energy}
E(u), E(v) \leq |||\kappa(\Hh)||| +c_L -c_R.
\end{equation} 

\subsection{Small cap data and a chain isomorphism in Morse homology}
\label{detect}

For an almost complex structure $J$ in $\JJ(M,\om)$, 
let $\hbar(J)$ be the infimum
over the symplectic areas of all nonconstant $J$-holomorphic 
spheres in $M$. We then set
$$
\hbar = \sup_{J \in \JJ(M,\om)} \hbar(J).
$$
This constant is positive and greater than or equal to $r(M,\om)$.

We now describe a chain isomorphism in Morse homology 
which can be constructed using cap data $\Hh$ that satisfies 
\begin{equation}\label{curvature inequality}
|||\kappa(\Hh)|||+c_L-c_R <\hbar.
\end{equation}
This chain map will be used in \S \ref{limit caps} to find 
central periodic orbits whose right Floer caps will, in turn, be used to detect the right 
asymptotic spanning disc of Proposition \ref{keycap} in \S \ref{finale}. 

Let $f$ be a Morse function on $M$ and $g$ a metric on $M$ such that the Morse complex, 
$(\CM(f), \p_g)$, is well-defined. Here $\CM(f)$ is the $\Z$-module generated by the 
critical points of $f$. The $\Z$-module generated by the elements of $\PP(H)$
is denoted by $\CF(H)$.

\begin{Proposition}\label{central caps}
Let $\Hh$ be  a generic choice of cap data for $(H,J)$.
If $|||\kappa(\Hh)|||+c_L-c_R <\hbar$, then there are two $\Z$-module homomorphisms
$$ \sl \colon \CM(f) \to \CF(H)$$ and $$\sr \colon \CF(H) \to \CM(f)$$ whose composition 
$$\Phi_{\Hh} = \sr \circ \sl \colon \CM(f)  \to \CM(f) $$
is a chain map which is chain homotopic to the identity.
\end{Proposition}

This result is strongly motivated by the work of Chekanov in \cite{ch}. The proof of Proposition \ref{central caps} is contained in \cite{ke2} where it appears 
as Proposition 2.4. While it is assumed there that 
$c_L=c_R=0$, the proof from \cite{ke2} extends easily to the present setting. The genericity assumption 
of this result concerns 
the almost complex structure $J$ as well as the families of almost complex structures appearing in the 
cap data $\Hh$. As usual, this assumption is included to ensure that the moduli spaces used to construct the maps are regular. These 
almost complex structures should also by chosen to lie in specific open sets of $\JJ(M, \om)$ 
so that inequality \eqref{curvature inequality} can be used to avoid bubbling.
These technical details, which are discussed in detail in \cite{ke2}, can be safely ignored in 
the present discussion.

Since the maps $\sr$, $\sl$ and $\Phi_{\Hh}$ play  important roles in the 
proof of Proposition \ref{keycap}, we will recall the relevant aspects of their constructions.
We begin with the map $\sl$. 

A left or right Floer cap will be called  {\bf short} if its energy is less than $\hbar$.
The subset of short elements in $\LL(x;\HH_L)$ will be denoted by $\LL'(x;\HH_L)$. 
Consider the space of left-half gradient trajectories;
$$
\ell(p) =\left\{ \alpha \colon (-\infty,0] \to M \mid \dot{\alpha} = - \nabla_g f(\alpha),\, \alpha(-\infty)=p \right\}.
$$
For a critical point $p$ of $f$ and an orbit $x$ in $\PP(H)$, set 
$$
\LL(p,x;f,\HH_L) =\{ (\alpha,u) \in \ell(p) \times \LL'(x;\HH_L) \mid \alpha(0)=u(-\infty)\}.
$$ 
For generic data, $\LL(p,x;f,\HH_L)$ is a smooth manifold and the local dimension of the 
component containing $(\alpha,u)$ is 
\begin{equation}
\label{dim}
\ind(p)-n - \CZ(x,u), 
\end{equation}
\cite{pss}. The  homomorphism
$
\sl \colon \CM(f) \to \CF(H)
$ 
is defined on each generator $p$ of $\CM(f)$ by
$$
\sl(p) = \sum_{x \in \PP(H)}  \# \LL_0(p, x ;f,\HH_L) x,
$$
where $\# \LL_0(p, x ;f,\HH_L)$ is the number of zero-dimensional components in 
$\LL(p, x ;f,\HH_L)$ counted with signs determined by a fixed coherent orientation. 
The shortness assumption implies that $\LL_0(p, x ;f,\HH_L)$  
is compact and so the map $\sl$ is well-defined. 

Next we consider the 
space of right-half gradient trajectories
$$
r(q) =\left\{ \beta \colon [0, +\infty) \to M \mid \dot{\beta} = - \nabla_g f(\beta),\, \beta(+\infty)=q \right\},
$$
and define
$$
\RR(x,q;\HH_R,f) =\{ (v, \beta) \in  \RR'(x;\HH_R)\times r(q) \mid v(+\infty) = \beta(0)\}.
$$
Here, $\RR'(x;\HH_R)$ is the set of short right Floer caps of $x$.
For generic data, each $\RR(x,q;\HH_R,f)$ is a smooth manifold, 
and the dimension of the 
component containing $(v, \beta)$ is $\CZ(x,\overleftarrow{v}) - \ind(q)+m$. 
Let $\RR_0(x,q ; \HH_L, f)$ be the set of zero-dimensional components in $\RR(x,q ; \HH_L, f)$.
The map $\Phi_{\Hh}$  is then defined by setting the coefficient 
of $q$ in $\Phi_{\Hh}(p)$ to be the integer
$$
 \sum_{x \in \PP(H)} \# \Big\{  ((\alpha, u), (v,\beta)) \in \LL_0(p,x;\HH_R,f)) \times \RR_0(x,q;\HH_R,f)) \mid [u\#v]=0 \Big\}.
$$ 

The map $\Phi_R \colon \CF(H) \to \CM(f)$ is determined by $\Phi_L$ and $\Phi_{\Hh}$
as follows. Let $V_L$ be the submodule of $\CF(H)$ generated by the orbits in $\PP(H)$ 
which appear in an element in the image of $\sl$ with a nonzero 
coefficient. The maps $\sl$ and $\Phi_{\Hh}$ uniquely determine 
the restriction of  $\sr$ to $V_L$. 
Setting $\sr=0$ on the complement of $V_L$ we obtain the 
full map. In particular, the coefficient of $q$
in $\sr(x)$ is the signed count of elements $(v,\beta) \in   \RR_0(x,q;\HH_R,f))$ for which there is an element 
$(u, \alpha)$ in some $\LL_0(p,x;\HH_R,f))$ such that $[u\#v]=0$.

\section{Proof of Proposition \ref{keycap}}

\subsection{Step 1: approximating $H_r$ by generic functions}

The results of the previous section can not be applied directly to $H_r$ because the elements 
of $\PP(H_r)$ are degenerate. To overcome this, 
we now approximate $H_r$ by a sequence of functions 
$H_k$ whose $1$-periodic orbits are nondegenerate. These functions are constructed 
explicitly in order to retain suitable control over their periodic orbits.

Let $F^0 \colon M \to \R$ be a Morse-Bott function 
with the following properties:
\begin{itemize}
\item The submanifold $N$ is a critical submanifold with index equal to $\codim(N)=k$.
\item All other critical submanifolds are isolated critical points of Morse index less than $\dim(M)=2m$.
\item On $U_r$,  $F^0=f^0(|p|)$ for some decreasing function $f^0 \colon [0,r] \to \R$.
\end{itemize}
Let $f_N \colon N \to \R$ be a Morse function with a unique local maximum at a point $Q$ in $N$.
Choose a bump function $\hat{\sigma} \colon [0,+\infty) \to \R$ such that $\hat{\sigma}(s)=1$ for $s$ near zero and 
$\hat{\sigma}(s)=0$ for $s \geq r/4$. Let $\sigma = \hat{\sigma}(|p|)$ be the corresponding function on $M$ with 
support in $U_{\scriptscriptstyle{r/4}}$ and set
$$
F= F^0 + \eps_N \cdot \sigma \cdot f_N.
$$
For a sufficiently small choice of $\eps_N>0$, $F$ is a Morse 
function whose critical points away from $U_{\scriptscriptstyle{r/4}}$ agree with those of 
$F^0$ and whose critical points in $U_{\scriptscriptstyle{r/4}}$ are precisely the critical points of $f_N$ on $N 
\subset M$.

Now let 
$$
H^0_k = H_r + \frac{1}{k} F.
$$
Each $H^0_k$ is also a Morse function with $\Crit(H^0_k) =\Crit(F)$. As well, $Q$ is the only critical point 
of $H^0_k$  with  Morse index equal to $2m$. For an interval $I \subset [0,R/2]$, we introduce the notation
$
U_{\scriptscriptstyle{I}} = \{ (q,p) \in U_{\scriptscriptstyle{R}} \mid |p| \in I \}.
$
When $k$ is sufficiently large, the $1$-periodic 
orbits of $H^0_k$ are either critical points or nonconstant orbits contained in
$U_{\scriptscriptstyle{(r/3, 2r/3)}}$. In fact, these nonconstant 
orbits lie in $U_{\scriptscriptstyle{[r/3+ \delta, 2r/3-\delta]}}$ 
for some $\delta>0$. This follows from the fact that
$dH^0_k$ converges to zero in the $C^{\infty}$-topology
along the boundary of $U_{\scriptscriptstyle{(r/3, 2r/3)}}$.

Perturbing each $H^0_k$ away from $\Crit(F)$, one obtains a 
sequence of  functions $H_k$  with the following properties
\begin{itemize}
\item $H_k \to H_r$ in $C^{\infty} (S^1 \times M)$.
\item The orbits  in $\PP(H_k)$ are nondegenerate and are of two types:  constant orbits 
which  coincide with the critical points of $F$,  and nonconstant orbits in $U_{\scriptscriptstyle{[r/3+ \delta, 2r/3-\delta]}}$ 
for some $\delta >  0$.
\item The constant periodic orbits equipped with their constant spanning discs have Conley-Zehnder indices
less than $m$ except for the constant orbit at the point $Q \in N$, which has Conley-Zehnder index equal to $m$.
\end{itemize}

The final detail to account for is the normalization condition. If we add the function 
$-\int^1_0 H_k(t, \cdot) \, \om^m$ to $H_k$, the resulting function is normalized and retains the 
properties described above. In particular, it determines the same Hamiltonian vector field. 
Abusing notation, this new normalized function will still be denoted by
$H_k$.

The following lemma provides a simple criteria for detecting nonconstant periodic orbits 
of $H_k$.
  
\begin{Lemma}\label{nonconstant}
If $x(t)$ is a $1$-periodic orbit of $H_k$ which admits a spanning disk $w$
such that $-\|H_k\|^- \leq \AC_{H_k}(x,w) < \|H_k\|^+$ and  $\CZ(x,w)=m$, then $x(t)$ is nonconstant.
\end{Lemma}

\begin{proof}
Arguing by contradiction, we assume that $x(t)=P$ for some point $P$ in $M$. 
The spanning disk $w$ then represents an element 
$[w]$ in $\pi_2(M)$, and we have 
\begin{equation*}\label{action}
\AC_{H_k}(x,w) = \int_0^1 H_k(t,P) \,dt - \om([w]).
\end{equation*}
Moreover, the point $P$ 
corresponds to a critical point of $F$ and  
$H_k$ is $C^2$-small near $P$, so our normalization of the Conley-Zehnder index yields
\begin{equation}\label{index}
\CZ(x,w) = \ind(P) -m-2c_1([w]).
\end{equation}

If $\om([w])=0$, then assumption \eqref{link} implies that $c_1([w]) \geq 0$.\footnote{This is the only point in the paper where we use assumption \eqref{link}.} It then follows from \eqref{index} that 
the Morse index of $P$ must be  $2m$. 
This implies that $P=Q$, since $Q$ is the unique fixed local maximum of $H_k$. However, the action of $Q$ 
with respect to a spanning disc $w$ with $\om([w])=0$ is equal to 
$\|H_k\|^+$. This is outside the assumed action range and hence a contradiction.

We must therefore have $\om([w])\neq 0$ and thus 
$$|\om([w])| \geq r(M,\om)> \|H_k\|.$$
For the case $\om([w])<0$, this implies that 
$$
\AC_{H_k}(x,w) \geq \int_0^1 H_k(t,P) \,dt + \|H_k\| \geq \|H_k\|^+
$$
which is a contradiction, as above. If  $\om([w])>0$, then
$$
\AC_{H_k}(x,w) \leq \int_0^1 H_k(t,P)\,dt  - \|H_k\|  = \int_0^1 H_k(t,P)\,dt - \|H_k\|^+ -\|H_k\|^-.
$$
So, either $\AC_{H_k}(P,w) < -\|H_k\|^-$ or $P=Q$. Both of these conclusions again contradict our hypotheses. Therefore  $x(t)$ must be nonconstant. 
\end{proof}

\subsection{Step 2: curve shortening}\label{curve shortening}

We now prove that the Hamiltonian path $\phi^t_{H_r}$
does not minimize the positive Hofer length in its homotopy class.
We also show that the same is true of the paths $\phi^t_{H_k}$ when $k$ is 
sufficiently large. 

For a Hamiltonian path $\psi_t$, let $[\psi_t]$ be
the class of Hamiltonian paths which are homotopic to $\psi_t$ relative to its endpoints. Denote the set of normalized functions 
which generate the paths in $[\psi_t]$ by 
$$
C^{\infty}_0([\psi_t]) =\{ H \in C^{\infty}_0(S^1 \times M) \mid [\phi^t_H \circ \psi_0] = [\psi_t]\}.
$$  
The Hofer semi-norm of $[\psi_t]$ is then defined by 
$$
\rho([\psi_t]) = \inf_{H \in C^{\infty}_0([\psi_t])} \{ \|H\| \}.
$$
The positive and negative Hofer semi-norms of $[\psi_t]$ are defined similarly as
$$
\rho^{\pm}([\psi_t]) = \inf_{H \in C^{\infty}_0([\psi_t])} \{ \|H\|^{\pm} \}.
$$
Clearly
$$
\rho([\psi_t]) \geq \rho^+([\psi_t]) + \rho^-([\psi_t]).
$$
In these terms,  the displacement energy of a subset $U \subset M$ is equal to
$$
e(U) = \inf_{\psi_t} \{ \rho([\psi_t]) \mid \psi_0=\id  \text{ and }\psi_1(U) \cap U = \emptyset \}.
$$

The following result is an easy application of Sikorav's curve shortening 
procedure. The proof follows very closely the proof of Proposition 2.1 
in \cite{schl}. 

\begin{Lemma}\label{shorten}
Let $H$ be an autonomous normalized Hamiltonian that is constant and equal to its minimal value on the 
complement of an open set $U \subset M$ which has finite displacement energy. 
If $\|H\|^+ > 2 e(U) $, then $$\|H\|^+ > \rho^+([\phi^t_H])  + \12 \|H\|^-.$$ 
In other words, $\phi^t_H$ does not minimize the positive Hofer semi-norm in its homotopy class.
\end{Lemma}

\begin{proof}
Let $\phi_t$ and $\psi_t$ be Hamiltonian paths and let $\varphi$ be a symplectomorphism.
The following properties of the positive and negative Hofer semi-norms are easily checked.
\begin{itemize}
  \item $\rho^{\pm}([\phi_t \circ \psi_t]) \leq \rho^{\pm}([\phi_t ] )+ \rho^{\pm}([ \psi_t])$ 
  \item $\rho^{\pm}([\phi_t \circ \psi]) = \rho^{\pm}([\phi_t])$
  \item $\rho^{\pm}([\psi^{-1} \circ \phi_t \circ \psi] )= \rho^{\pm}([\phi_t])$
  \item $\rho^+([\phi^{-1}_t]) = \rho^-([\phi_t]).$
\end{itemize}

Now choose a Hamiltonian path $\psi_t$ starting at the identity such that 
$$
\psi_1(U) \cap U = \emptyset.
$$ 
The path $\phi^t_H$ can then be factored as follows.
\begin{eqnarray*}
\phi^t_H & = & \left( \phi^{t/2}_H \circ \psi_t \circ \phi^{t/2}_H \circ \psi_t^{-1}\right) \circ \left(\psi_t \circ (\phi^{t/2}_H)^{-1} \circ \psi_t^{-1} \circ \phi^{t/2}_H \right)\\
{} & = & b_t \circ a_t.
\end{eqnarray*}
Hence,
$$
\rho^+([\phi^t_H])  \leq \rho^+([a_t]) + \rho^+([b_t]).
$$
For the first summand on the right, we have
\begin{eqnarray*}
 \rho^+([a_t])& = & \rho^+([\psi_t \circ (\phi^{t/2}_H)^{-1} \circ \psi_t^{-1} \circ \phi^{t/2}_H]) \\
{} & \leq&  \rho^+([\psi_t] )+ \rho^+([ (\phi^{t/2}_H)^{-1} \circ \psi_t^{-1} \circ \phi^{t/2}_H] )\\
{} & =&  \rho^+([\psi_t] )+ \rho^+([ (\phi^{1/2}_H)^{-1} \circ \psi_t^{-1} \circ \phi^{1/2}_H] )\\
{} & =&  \rho^+([\psi_t])  + \rho^+([ \psi_t^{-1}] )\\
{} & = &  \rho^+([\psi_t])  + \rho^-([ \psi_t] )\\
{} &\leq&  \rho([\psi_t]),
\end{eqnarray*}
and for the second summand we have
\begin{eqnarray*}
\rho^+([b_t]) & = & \rho^+([\phi^{t/2}_H \circ \psi_t \circ \phi^{t/2}_H \circ \psi_t^{-1}])\\
{} & = &  \rho^+([\phi^{t/2}_H \circ \psi_1 \circ \phi^{t/2}_H \circ \psi_1^{-1}])\\
{} & \leq &  \left\| \12 H + \12 H \circ \psi_1^{-1} \circ (\phi^{t/2}_H)^{-1}\right\|^+\\
{} & = &  \left\| \12 H \circ \phi^{t/2}_H + \12 H \circ \psi_1^{-1} \right\|^+\\
{} &=& \max_{p \in M} \left(\12 H(p) + \12 H \circ \psi_1^{-1}(p)\right)\\
{} &=& \12 \|H\|^+ - \12 \|H\|^- .
\end{eqnarray*}
Together, these inequalities imply that for any Hamiltonian path $\psi_t$ which displaces $U$
we have  
$$
\rho^+([\phi^t_H]) \leq  \rho([\psi_t] )+\12 \|H\|^+ - \12 \|H\|^-.
$$
Taking the infimum over all such paths we get 
\begin{eqnarray*}
  \rho^+([\phi^t_H]) &\leq&  e(U) +\12 \|H\|^+ - \12 \|H\|^- \\
  {} &<& \|H\|^+   - \12 \|H\|^-.
\end{eqnarray*}
\end{proof}

By construction, we have  $\|H_r\|^+ =A  > 2 e(U_r)$. Together with Lemma \ref{shorten} this implies that 
$\phi^t_{H_r}$ does not minimize the positive Hofer length in its homotopy class. In other words, 
there is a function $G$ in $C^{\infty}_0([\phi^t_{H_r}])$ such that 
\begin{equation}\label{g2h}
\|H_r\|^+ = \|G\|^+ + 2\eps
\end{equation}
for some $\eps>0.$ We now show that for sufficiently large $k$, the paths $\phi^t_{H_k}$ can 
be shortened in their respective 
homotopy classes by at least $\epsilon$.

To see this, consider the Hamiltonian path $\phi^t_{H_k} \circ (\phi^t_{H_r})^{-1}$ which is generated 
by the function 
$$
F_k =  H_k - H_r\circ \phi^t_{H_r} \circ (\phi^t_{H_k})^{-1}.
$$
These functions clearly converge to zero in the $C^{\infty}$-topology. 
The path 
$$
\phi^t_{H_k} \circ (\phi^t_{H_r})^{-1} \circ \phi^t_G
$$
is homotopic to $\phi^t_{H_k}$ and is generated by the function
$$
G_k = F_k +G \circ (\phi^t_{F_k})^{-1}.
$$
Hence, we have functions  $G_k$  in $C^{\infty}_0([\phi^t_{H_k}])$ such that  
$                
\|G_k\|^+ = \|F_k\|^+ + \|G\|^+. 
$ 
For large enough $k$ we then have
 \begin{equation}
\label{k G H}                  
\|G_k\|^+\leq \|H_k\|^+ - \eps. 
\end{equation}  

\subsection{Step 3: nontrivial linear right Floer caps}
\label{limit caps}

Fix a family of almost complex structures $J_k$  for each $H_k$ such that the $J_k$ converge to 
$J$ in $\JJ_{S^1}(M,\om)$.
As in Example \ref{ex:linear}, set  
$$
\overline{\HH}_k= \big((\eta(s) -1)\|H_k\|^- + \eta(s) H_k,0, J_{k,s}\big),
$$ 
where the $(\R \times S^1)$-families of almost complex structures $J_{k,s}$  converge to a compact  homotopy
$J_s$ from some $J^-$ to $J$.  The linear homotopy triples $\overline{\HH}_k$ then 
converge to the linear homotopy triple 
$\overline{\HH}_r= (\overline{H}_s, 0, J_s)$ for $(H_r,J)$. 

\begin{Proposition}\label{rlc}
For large enough $k$, there is a  nonconstant  $1$-periodic orbit $x_k$ of $H_k$ and a right 
Floer cap $v_k$  in $\RR(x_k;  \overline{\HH}_k)$ such that 
\begin{equation}\label{energy}
E(v_k) <\|H_k\|< r(M,\om)
\end{equation}
and  
\begin{equation}\label{action}
-\|H_k\|^- \leq \AC_{H_k}(x_k, \overleftarrow{v_k})< \|G_k\|^+ \leq \|H_k\|^+ - \eps.
\end{equation}
\end{Proposition}

\begin{proof}

As shown above, for large enough $k$ there is a function
$G_k$ such that the Hamiltonian path $\phi^t_{G_k}$ is homotopic to 
$\phi^t_{H_k}$, relative endpoints, and $\|G_k\|^+ \leq \|H_k\|^+ +\epsilon.$
Applying Proposition \ref{good caps} to $G_k$, we get a $J_k$ in $\JJ_{S^1}(M, \om)$ 
and a homotopy triple $\HH_{G_k}$ for $(H_k,J_k)$ such that 
\begin{equation}
\label{k G}
|||\kappa(\HH_{G_k})|||^+ + c_{L,k} \leq \|G_k\|^+.
\end{equation}

We now consider the following cap data for $(H_k,J_k)$,
$$\Hh_k = (\HH_{G_k}, \overline{\HH}_k).$$ 
By inequalities \eqref{k G H}, \eqref{k G}, and the curvature norm estimates for linear homotopy 
triples derived in Example \ref{ex:linear}, we have
\begin{equation}
\label{k inequality}
|||\kappa(\HH_{G_k})|||^+ +|||\kappa(\overline{\HH}_k)|||^-+c_{L,k}  - c_{R,k}  \leq \|G_k\|^+ + \|H_k\|^-.
\end{equation}
By construction (see inequalities \eqref{assumption r} and \eqref{Hr bounds}) we also have 
$\|H_r\|< r(M,\om)$. Hence, for sufficiently large $k$, \eqref{k inequality} implies that  
\begin{equation}
\label{k key}
|||\kappa(\HH_{G_k})|||^+ +|||\kappa(\overline{\HH}_k)|||^-  +c_{L,k} - c_{R,k}  < \|H_k\|< r(M,\om).
\end{equation}
From this point on, we will assume that $k$ is large enough for this inequality to hold.

Since $r(M, \om) \leq \hbar$, inequality \eqref{k key} allows us to apply Proposition \ref{central caps} 
to the homotopy data $\Hh_k$. In particular,  for any 
Morse-Smale pair $(f,g)$ on $M$ we can construct two $\Z$-module homomorphisms
$$ \srl \colon \CM(f) \to \CF(H_k)$$ and $$\sll \colon \CF(H_k) \to \CM(f)$$ whose composition is a chain map
$$\Phi_{\Hh_k} \colon \CM(f) \to \CM(f)$$
which is chain homotopic to the identity.

For simplicity, we choose the Morse-Smale pair $(f,g)$ so that the function $f$  
has a unique local (and hence global) maximum at a point $q \in M$.
Standard arguments imply that $q$ is the unique nonexact cycle of degree $2m$ in the 
Morse complex $(\CM(f), \p_{g})$, and so
$$\Phi_{\Hh_k}(q) =q.$$

Let $V_{L,k}$ be the submodule of $\CF(H_k)$ generated by $1$-periodic orbits of $H_k$ which appear 
in an element of the image of $\srl$ with a nonzero coefficient. Let $K_{R,k}$ be the submodule of $\CF(H_k)$ generated by periodic orbits which lie in the kernel 
of $\sll$ and let $p_k \colon V_{L,k} \to V_{L,k} / K_{R,k}$ be the projection map. We then have
\begin{equation*}
\label{break}
\Phi_{\Hh_k} = \sll \circ p_k \circ \srl.
\end{equation*}
It follows from the definitions of these maps that any periodic orbit which appears in the image of $p_k\circ \srl$ is central
with respect to $\Hh_k$.

Let $$X_k = p_k \circ \srl(q).$$  By the construction of $\Phi_{\Hh_k}$, $X_k$ is a finite sum of the form $$X_k=\sum n^j_k x^j_k$$
where the $n^j_k$ are nonzero integers and the $x^j_k$ are 
central $1$-periodic orbits of $H_k$.

Since $X_k$
gets mapped to $q$ under $\srl$, the moduli space  
$$
\RR_0(X_k,q;\overline{\HH}_k,f) = \bigcup_j\RR_0(x^j_k,q;\overline{\HH}_k,f),
$$
which determines the image $\sll(X_k)$, must be nonempty.  
Choose a 
$(v_k, \sigma_k)$ in $\RR_0(X_k,q;\overline{\HH}_k,f)$  for each $k$. The caps $v_k$ belongs to $\RR(x_k;\overline{\HH}_k)$  for some $x_k$ in $\PP(H_k)$ which appears in  $X_k$ with a nonzero coefficient. 
Moreover, $v_k$ is part of a central pair for $x_k$ with respect to $\Hh_k$, and so by \eqref{action-bounds}, \eqref{k G H} and \eqref{k G}, we have 
\begin{equation*}\label{}
-\|H_k\|^- \leq \AC_{H_k}(x_k,\overleftarrow{v_k})\leq \|H_k\|^+-  2\eps.
\end{equation*}
Inequality  \eqref{uniform-energy} together with \eqref{k key} yields the desired  uniform energy bound
\begin{equation*}\label{}
E(v_k) \leq |||\kappa(\Hh_k)||| +c_{L,k}-c_{R,k}  <  \|H_k\| < r(M,\om).
\end{equation*}

It only remains to show that the orbits $x_k$ are nonconstant. Each 
$x_k$ appears in $p_k \circ \srl(q)$ with a nonzero coefficient.  Hence, there is a pair of maps $(\alpha_k, u_k)$ in
$\LL_{[0]}(q,x_k;f,\HH_{G_k})$ such that $u_k$ is part of a central pair for $x_k$
with respect to $\Hh_k$. The existence of the regular pair $(\alpha_k,u_k)$
together with the  
dimension formula  for $\LL_{[0]}(q,x_k;f,\HH_{G_k})$, \eqref{dim}, implies 
that $$\CZ(x_k, u_k)= \ind(q)-m =m.$$ 
Since $u_k$ is part of a central pair for $x_k$
the  action
$\AC_{H_k}(x_k,u_k)$ satisfies the same bounds, \eqref{action}, as  $\AC_{H_k}(x_k,v_k)$, i.e.,  
\begin{equation*}
\label{ }
-\|H_k\|^- \leq \AC_{H_k}(x_k,u_k)\leq \|H_k\|^+-  2\eps.
\end{equation*}
Lemma \ref{nonconstant} then implies that the orbits $x_k$ are nonconstant
and the proof of Proposition \ref{rlc} is complete.

\end{proof}

\subsection{Step 4: A nonconstant limit of linear right Floer caps}\label{finale}

Let $\CC$ be the closed subset of $C^{\infty}(\R \times S^1, M)$ consisting of maps 
$v \colon \R \times S^1 \to M$ such that $v(0,t)$ is a contractible loop in $M$. We consider this space 
as being equipped with the $C^{\infty}_{loc}$-topology.

By Proposition \ref{rlc} we have a sequence of nonconstant periodic orbits 
$x_k \in \PP(H_k)$ and a sequence of right Floer caps $v_k \in \RR(x_k, \overline{\HH}_k)$
which satisfy \eqref{energy} and \eqref{action}.  
The linear homotopy triples $\overline{\HH}_k$ were chosen so that they converge to $\overline{\HH}_r= (\overline{H}_s,0,J_s)$. 
Together with uniform energy bound \eqref{energy}, this implies that there is a subsequence of the $v_k$ 
which converges in $\CC$ to a map $\tilde{v}$. This limiting map $\tilde{v}$ 
is a solution of the equation
\begin{equation}\label{limit}
\partial_s \tilde{v}  + \overleftarrow{J_s}(\tilde{v})(\partial_t\tilde{v} - X_{\overleftarrow{\overline{H}}_s}(\tilde{v}))=0,
\end{equation}
for $$\overleftarrow{\overline{H}}_s = (\eta(-s)-1)B + \eta(-s)H_r.$$
It also satisfies 
\begin{equation}
\label{energy limit}
0 \leq E(\tilde{v})< r(M,\om).
\end{equation}

The map $\tilde{v}$ may or not be constant.
To find the the periodic orbit and the right asymptotic spanning disc 
of Proposition \ref{keycap}, we need to consider both possibilities.

\subsubsection{Case 1: a nonconstant limit}\label{case 1}

We assume here that the subsequence, which we still denote by $v_k$, converges to 
a nonconstant solution $\tilde{v}$ of \eqref{limit}. In this case, the map $\tilde{v}$ will be
the asymptotic right spanning disc of Proposition \ref{keycap} and we will write $\tilde{v}=v$.

The energy bound  \eqref{energy limit} implies that the limit $v(+\infty)=\lim_{s \to +\infty}v(s,t)$ is a point in $M$.
It also implies that there is a sequence $s_j^- \to -\infty$ such that  $v(s_j^-,t)$ converges 
to some $y(t)$ in $\PP(H_r).$  For simplicity we assume that the sequence $s_j^-$
is monotone decreasing and that $s_1^- < -1$.

It remains for us to show that the limiting periodic orbit 
$y$ is nonconstant and that there is an $\epsilon>0$ such that \eqref{key1} holds for all $j$, i.e., 
\begin{equation*}\label{jae}
E^j(v) -B \leq \AC^j_{H_r}(v) \leq A -\epsilon. 
\end{equation*}
We begin by proving that  \eqref{key1} holds for $\epsilon = \frac{1}{2} (\|H_r\|^+ - \|G\|^+)$. 
By Remark \ref{new terms}, we have $E^j(v) = \overleftarrow{E}^{-s_j^-}_{s_j^-}(v)$ and 
$\AC^j_{H_r}(v) = \overleftarrow{\AC}^{-s_j^-}_{s_j^-}(v)$.
Since the $v_k$ converge to $v$ in the $C^{\infty}_{loc}$-topology, and $H_k$ converges to
$H_r$ in the $C^{\infty}$-topology, it suffices to show that for large enough $k$ we have 
\begin{equation}\label{jk}
\overleftarrow{E}^{-s_j^-}_{s_j^-}(v_k)-\|H_k\|^- \leq \overleftarrow{\AC}^{-s_j^-}_{s_j^-}(v_k) \leq \|H_k\|^+-\eps.
\end{equation}

The first of these inequalities follows immediately from equation \eqref{action energy right}.
In particular, this identity implies that
\begin{eqnarray*}
\overleftarrow{\AC}^{-s_j^-}_{s_j^-}(v_k) & \geq &  \overleftarrow{E}^{-s_j^-}_{s_j^-}(v_k) -\|H_k\|^-
+\int_0^1 \int^{-s_j^-}_{s_j^-} \kappa(\overline{\HH}_k)(s,t,\overleftarrow{v_k}(s,t))\,ds\,dt\\
{} & \geq & \overleftarrow{E}^{-s_j^-}_{s_j^-}(v_k)-\|H_k\|^-,  
\end{eqnarray*}
since $ \kappa(\overline{\HH}_k) = \dot{\eta}(s)(\|H_k\|^- + H_k) \geq 0$. 

To prove the second inequality in \eqref{jk}, we first note that for $i>j$ inequality \eqref{action difference}
yields 
\begin{equation*}
\overleftarrow{\AC}^{-s_i^-}_{s_i^-}(v_k) - \overleftarrow{\AC}^{-s_j^-}_{s_j^-}(v_k) \geq
\int_0^1 \int_{[s_i^-,-s_j^-] \cup [-s_j^-, -s_i^-]} \kappa(\overline{\HH}_k)(s,t,\overleftarrow{v_k}(s,t))\,ds\,dt.
\end{equation*}
Hence, for each $k$, the sequence $\overleftarrow{\AC}^{-s_j^-}_{s_j^-}(v_k)$ is nondecreasing and 
$$
\AC_{H_k}(x_k,\overleftarrow{v_k})= \lim_{i \to \infty}\overleftarrow{\AC}^{-s_i}_{s_i}(v_k) \geq \overleftarrow{\AC}^{-s_j^-}_{s_j^-}(v_k).
$$
By \eqref{action}, this yields  
$$
\overleftarrow{\AC}^{-s_j^-}_{s_j^-}(v_k) \leq \|H_k\|^+ -\eps.
$$
Thus, \eqref{jk} holds and we have established inequality \eqref{key1}.

Finally we must show that the periodic orbit $y$ is nonconstant. This is easily derived from \eqref{key1}
as follows.
Set 
$$y^{\scriptscriptstyle{[s]}}(t) = v(s,t)$$ and consider the annulus 
$$
v^{\scriptscriptstyle{[s]}} = v|_{\scriptscriptstyle{[-s,s] \times S^1}}.
$$
Since  $y^{\scriptscriptstyle{[s_j^-]}} \to y$ and $y^{\scriptscriptstyle{[-s_j^-]}} \to p$, 
for large values of $j$ the 
annuli $v^{\scriptscriptstyle{[s_j^-]}}$ can be extended and reparameterized to form spanning discs
for $y \in \PP(H)$ in a fixed homotopy class. These extensions 
can be made arbitrarily small for sufficiently large $j$. Hence, by inequality \eqref{key1} and the 
assumption that $v$ is nonconstant,
we can choose such a spanning disc $w$ for $y$ such that 
\begin{equation}\label{interval}
-B < \AC_{H_r}(y,w) < A.
\end{equation}
Assume now that $y$ is a constant periodic orbit, i.e., $y(t)=P$ for some critical point 
$P$ of $H_r$. 
Then $w$ represents a class $[w] \in \pi_2(M)$ and 
\begin{equation}\label{false action}
\AC_{H_r}(y,w)= H_r(P)-\om([w]).
\end{equation}

If $\om([w])=0$, then \eqref{interval} and \eqref{false action} imply that $P$ must be a critical point of $\H_r$ with critical value in 
$(-B,A)$. Since there is no such critical point, we must have $\om([w]) \neq 0$ and hence $$|\om([w])|\geq r(M,\om) > \|H_r\|=A+B.$$
However, this implies that $\AC_{H_r}(y,w)$ fails to lie in the interval $(-B,A)$, which contradicts 
\eqref{interval}. The orbit $y$ must therefore be nonconstant.

\subsubsection{Case 2: a constant limit}

We now assume that the maps $v_k$ converge in $\CC$ to a constant map 
$\tilde{v}(s,t) = \widetilde{P}$. In this case, we will adapt a topological argument from \cite{gi:coisotropic}
to prove that there is a sequence 
$\tau_k \to -\infty$, such that $v_k(s + \tau_k,t)$ converges to a nonconstant 
solution $v$ of the equation 
\begin{equation}
\label{floer}
\p_s v + J(v)(\p_t v -X_{H_r}(v))=0.
\end{equation}
This will be the right asymptotic spanning disc of Proposition \ref{keycap}.

To detect this map, 
we first pass to a subsequence of the $v_k$ whose negative asymptotic limits converge
to a nonconstant element of $\PP(H_r)$.
Recall that $x_k = v_k(-\infty)$ is a nonconstant 
$1$-periodic of $H_k$.
Since the $x_k$ are nonconstant, they are contained in the 
region $U_{\scriptscriptstyle{[r/3+ \delta ,2r/3 - \delta]}}$. 
By Arzela-Ascoli, there is a convergent subsequence of the $x_k$ 
that converges to some $x \in \PP(H_r)$. Since it 
is contained in $U_{\scriptscriptstyle{[r/3+ \delta ,2r/3 - \delta]}}$, the orbit $x$ is also nonconstant. 
From this point on, we restrict
our attention to a subsequence of the $v_k$ for which the $x_k$ converge to $x$. 
For simplicity, this subsequence will still be denoted 
by $v_k$.

There is a natural action of $\R$ on $\CC$ defined by
$\tau \cdot v(s,t) = v(s + \tau,t).$ 
We set
$$\Gamma(v_k)=\{\tau \cdot v_k \mid \tau \in \R\},$$
and define $\Sigma$ to be the set of 
limits of all convergent sequences of the form  
$$v = \lim_{k \to \infty}\tau_k \cdot v_k.$$

There are two continuous maps on $\Sigma$ which will be useful in what follows. The first is 
the evaluation map $\ev \colon \Sigma \to M$ defined by 
$$
\ev(v) = v(0,0).
$$
The second map is the function $\overleftarrow{\AC}^0_{-\infty}\colon  \Sigma \to \R$
which is defined, as in \S \ref{reverse action}, by
\begin{equation*}
 \overleftarrow{\AC}^0_{-\infty}(v) = \int_0^1 H_r(\overleftarrow{v}(0,t)) \, dt - \int_0^1 \int_{-\infty}^0 \om(\p_s \overleftarrow{v}, \p_t \overleftarrow{v}) \,ds\,dt.
\end{equation*}

\begin{Lemma}\label{action bounds}
For every $v$ in $\Sigma$, $$-B \leq \overleftarrow{\AC}^0_{-\infty}(v) \leq A-\epsilon.$$
\end{Lemma}
\begin{proof}
By the definition of $\Sigma$, we have $v =\lim_{k\to \infty} \tau_k \cdot v_k$ for some sequence $\tau_k$.
Hence, 
$$\overleftarrow{\AC}^0_{-\infty}(v) = \lim_{k \to \infty} \overleftarrow{\AC}^0_{-\infty}(\tau_k \cdot v_k)$$
and it suffices to show that for sufficiently large $k$ we have 
\begin{equation}
\label{k estimate}
\overleftarrow{E}^0_{-\infty}(\tau_k \cdot v_k)-\|H_k\|^- \leq \overleftarrow{\AC}^0_{-\infty}(\tau_k \cdot v_k)\leq \|H_k\|^+ -\epsilon.
\end{equation}
The proof of these inequalities is entirely similar to the  proof
of  \eqref{jk}. 
In particular, \eqref{action energy right} implies that
\begin{equation*}
\overleftarrow{\AC}^0_{-\infty}(\tau_k \cdot v_k) \geq  \overleftarrow{E}^0_{-\infty}(\tau_k \cdot v_k)-\|H_k\|^-.  
\end{equation*}
On the other hand, \eqref{action difference} yields
\begin{equation*}
\AC_{H_k}(x_k,\overleftarrow{\tau_k \cdot v_k}) - \overleftarrow{\AC}^0_{-\infty}(\tau_k \cdot v_k)
\geq
\int_0^1 \int_0^{+\infty} \kappa(\overline{\HH}_k)(s,t,\overleftarrow{\tau_k \cdot v_k}) \,ds\, dt,
\end{equation*}
and by \eqref{action}, we then have   
$$
\|H_k\|^+ -\eps \geq \AC_{H_k}(x_k,\overleftarrow{v_k})= \AC_{H_k}(x_k,\overleftarrow{\tau_k \cdot v_k}) \geq \overleftarrow{\AC}^0_{-\infty}(\tau_k \cdot v_k).
$$
\end{proof}

\begin{Lemma}\label{belong}
Every element of $\Sigma$ is a solution of  
\eqref{floer} with energy less that $r(M,\om)$. 
\end{Lemma}
\begin{proof}

Let $v =\lim_{k\to \infty} \tau_k \cdot v_k$.
By \eqref{energy}, the energy of each $v_k$ is less than $r(M, \om)$. 
Since  $E(v_k)= E(\tau_k \cdot v_k)$, 
the  energy of $v$ is also less than $r(M, \om)$.
It only remains to show that $v$ is a 
solution of \eqref{floer}.

Recall that, $\overleftarrow{\overline{H}}_s$ is a compact homotopy from 
$H_r$ to $-B$.  If $\tau_k \to -\infty$, 
then $v$  is clearly a solution of 
\eqref{floer}. If the sequence of shifts $\tau_k$ is bounded, then $v$ is equal to the constant map 
$\tilde{v}(s,t)= \widetilde{P}$. Since $\tilde{v}$ is also a solution
 of \eqref{limit}, 
we must have $X_{\overleftarrow{\overline{H}}_s}(\widetilde{P})=0$ for all $s \in \R$. In other words, $\widetilde{P}$ is a critical point of $H_r$ and hence must lie in $U_{\scriptscriptstyle{r/3}} \cup U_{\scriptscriptstyle{[2r/3,+\infty)}},$ where $U_{\scriptscriptstyle{[2r/3, +\infty)}}$ denotes the complement of $U_{\scriptscriptstyle{2r/3}}$ in $M$.
Lemma \ref{action bounds} implies that 
\begin{equation*}
-B \leq  \overleftarrow{\AC}^0_{-\infty}(\tilde{v}) = H_r(\widetilde{P}) \leq A-\epsilon.
\end{equation*} 
Hence, $\widetilde{P}$ belongs to $U_{\scriptscriptstyle{[2r/3,+\infty)}}$. 
On this set $\overleftarrow{\overline{H}}_s = H_r =-B$, and so 
$v$ is a trivial solution of \eqref{floer}.

Finally, when the shifts $\tau_k \to \infty$, the limit $v$ is a solution of 
\begin{equation*}
\p_s v + J^-(v)\p_t v=0
\end{equation*}
with energy less than $r(M,\om)$. Any such map can be uniquely extended to 
a holomorphic sphere with the same energy. Since $r(M,\om)<\hbar$, the almost complex structure $J^-$ can be chosen, 
at the outset, to satisfy
$\hbar(J^-)>\|H_r\|$. The map $v$ must therefore be constant. 
Lemma \ref{action bounds} implies that the constant maps in $\Sigma$ all lie in 
$U_{\scriptscriptstyle{[2r/3, +\infty)}}$. 
Hence, $v$ is again a trivial solution of \eqref{floer}.

\end{proof}

\begin{Lemma}
\label{nonincrease}
The function $\tau \mapsto \overleftarrow{\AC}^0_{-\infty} (\tau \cdot v)$
is nonincreasing. It is strictly decreasing unless $v$ belongs to $\PP(H_r)$.
\end{Lemma}
Here, the elements of $\PP(H_r)$ are
identified with elements of $\CC$ that do not depend on $s$. 
\begin{proof}
For $\tau'>\tau$, a simple computation yields
\begin{equation}
\label{a0}
\overleftarrow{\AC}^0_{-\infty} (\tau' \cdot v) - \overleftarrow{\AC}^0_{-\infty} (\tau \cdot v) = -\int_0^1 \int_{-\tau'}^{-\tau} \om(\p_s \overleftarrow{v}, J(\overleftarrow{v}) \p_s \overleftarrow{v}) \,ds\,dt.
\end{equation}
Since the integrand is nonnegative the function $\overleftarrow{\AC}^0_{-\infty} (\tau \cdot v)$
is nonincreasing. 

If  $\overleftarrow{\AC}^0_{-\infty} (\tau' \cdot v)  = \overleftarrow{\AC}^0_{-\infty} (\tau \cdot v)$ for $\tau'>\tau$, then \eqref{a0} implies
that $\p_s v = 0$ for $s \in (\tau, \tau')$. By Lemma \ref{belong}, $v$ is a solution of \eqref{floer}, and hence $v(s,t)=v(t)$ is
a $1$-periodic orbit of $H_r$ for $s \in (\tau, \tau')$. The Unique Continuation Theorem
of \cite{fhs} then implies that $v(s,t)=v(t)$ for all values of $s$.
\end{proof}
Following \cite{gi:coisotropic} we now prove: 
\begin{Lemma}\label{properties} 
The set $\Sigma$ has the following properties.
\begin{enumerate}
\renewcommand{\theenumi}{\roman{enumi}}
\renewcommand{\labelenumi}{(\theenumi)} 
  \item the point $\widetilde{P}$ and the nonconstant $1$-periodic orbit $x(t)$ belong to $\Sigma$;
  \item the subsets  $\Gamma(v_k) \subset \CC$ converge to $\Sigma$ in the Hausdorff topology;
  \item the set $\Sigma$ is connected, compact and preserved by the $\R$-action on $\CC$;
   \item The action of $\R$ on  $\Sigma$ is nontrivial.
 \end{enumerate}
\end{Lemma}

\begin{proof}
The first two properties follow almost immediately from the definition of 
$\Sigma$. The same is true of the fact that $\Sigma$ is invariant under the $\R$-action.

The compactness of $\Sigma$ follows from Lemma \ref{belong} and the fact that $\Sigma$ is closed.
In particular, the subset of $\CC$ consisting 
of solutions of \eqref{floer} with energy less than $\hbar$ is itself compact
by the usual Floer compactness theorem.

To prove that $\Sigma$ is connected, consider any two disjoint 
open sets in $\CC$, $\UU_1$ and $\UU_2$, which cover $\Sigma$. Let $\Sigma_{\widetilde{P}}$ be the 
component of $\Sigma$ which contains $\widetilde{P}$ and suppose that $\Sigma_{\widetilde{P}}  \subset \UU_1$. 
By (ii), the $\Gamma(v_k)$ are contained in $\UU_1 \cup  \UU_2$ for all sufficiently large $k$. Since 
the $\Gamma(v_k)$ are connected and $\widetilde{P}$ is a limit point of the $\Gamma(v_k)$,
they must be contained in $\UU_1$ for large $k$. Thus, $\Sigma\cap \UU_2=\emptyset$
and it follows that $\Sigma$ must be connected. 

To prove (iv), we note that (iii) implies that $ev(\Sigma)$ is connected. Since $\widetilde{P}$ belongs
to $U_{\scriptscriptstyle{[2r/3, +\infty)}}$ and $x(t)$ belongs to $U_{\scriptscriptstyle{[r/3+\delta ,\, 2r/3 - \delta]}}$ there must be some 
$v$ in $\Sigma$ for which  $ev(v)$ belongs to $U_{\scriptscriptstyle{(2r/3 - \delta, 2r/3)}}$. There are no $1$-periodic 
orbits on the level sets in $U_{\scriptscriptstyle{(2r/3 - \delta, 2r/3)}}$, so the loop $v(0,t)$ is not in $\PP(H_r)$. Hence, $v$ is not a fixed point of the $\R$-action by Lemma \ref{nonincrease}.

\end{proof}

We now consider the set 
$$\Sigma_{\scriptscriptstyle{\min}} = \{ v \in \Sigma \mid \overleftarrow{\AC}^0_{-\infty}(v) =-B\}$$ 
The properties of $\Sigma$ established above imply that $\Sigma_{\scriptscriptstyle{\min}}$ is comprised of elements in $\PP(H_r)$.
In particular, for $v \in \Sigma_{\scriptscriptstyle{\min}}$  choose a $\tau<0$. Lemmas 
\ref{action bounds}, \ref{nonincrease}, and \ref{properties} yield
$$
-B = \overleftarrow{\AC}^0_{-\infty}(v) \leq \overleftarrow{\AC}^0_{-\infty}(\tau \cdot v) \geq -B.
$$
The second  statement of Lemma \ref{nonincrease} then implies that $v$ belongs $\PP(H_r)$.

Note that the constant elements of $\Sigma_{\scriptscriptstyle{\min}}$ 
take values in the set $U_{\scriptscriptstyle{[2r/3, +\infty)}}$. Identifying $U_{\scriptscriptstyle{[2r/3, +\infty)}}$
with the space of constant maps in $\CC$ which take values in $U_{\scriptscriptstyle{[2r/3, +\infty)}}$,
we define $$C = U_{\scriptscriptstyle{[2r/3, +\infty)}} \cap \Sigma_{\scriptscriptstyle{\min}} =  U_{\scriptscriptstyle{[2r/3, +\infty)}} \cap \Sigma.$$ The set $C$ is a compact subset of $\Sigma$. By property
(iv) of Lemma \ref{properties}, $C$ is also a proper subset of $\Sigma$. Most importantly, $C$ is nonempty 
because it contains $\widetilde{P}$.

\begin{Lemma}
\label{connected}
The set $C$ is a union of connected components of $\Sigma_{\scriptscriptstyle{\min}}$.
\end{Lemma}

\begin{proof}
If one assumes the contrary, then there is a sequence of nonconstant periodic orbits $x^-_k  \in \Sigma_{\scriptscriptstyle{\min}} \ssminus C$
which converges to a point of $C$. This is a contradiction since the  nonconstant orbits are contained in the closure of $U_{\scriptscriptstyle{2r/3-\delta}}$.
\end{proof}

Let $\NN_c = \{ w \in \Sigma \mid \overleftarrow{\AC}^0_{-\infty}(v) < -B + c\}$.
Fix a connected component $C_0$ of $C$  and let $\VV_c$ be the component 
of $\NN_c$ which contains $C_0$. 

\begin{Lemma}
\label{basis}
For any open set $\VV$ in $\Sigma$  which contains $C_0$, there is a $c>0$ such that  $\VV_c \subset \VV.$ 
\end{Lemma}

\begin{proof}
Assume the contrary. Then there is neighborhood $\VV \supset C_0$ and a sequence $c_i \to 0^+$ such that
$\VV_{c_i}$ is not contained in $\VV$.  Let $v_i$ be an element in $\VV_{c_i} \ssminus \VV$. Since $\Sigma$ 
is compact, there is a subsequence of the $v_i$ which converges to an element of $\bigcap_i \overline{\VV}_{c_i} \ssminus \VV$.
On the other hand,  $\bigcap_i \overline{\VV}_{c_i}$ is a connected subset of 
$\Sigma_{\scriptscriptstyle{\min}}$ which contains $C_0$. This contradicts the fact that $C_0$ is a connected
component of $\Sigma_{\scriptscriptstyle{\min}}$.
  
\end{proof}

We can now complete the proof of Proposition \ref{keycap} in the present case. 
By Lemma \ref{basis}, we can find a constant $c>0$ such that $\VV_c \subset ev^{-1}(U_{\scriptscriptstyle{(2r/3-\delta, +\infty)}}).$ Either $\VV_c \cap C = C_0$ or $\VV_c \cap C$ is disconnected. In both cases, the fact that
 $\Sigma$ is connected and contains the nonconstant orbit $x(t)$ implies that $$\VV_c \ssminus C \neq \emptyset.$$ 
Let $v$ be any map in $\VV_c \ssminus C$. We will show that $v$ is a right asymptotic spanning disc with the 
desired properties.

Properties \eqref{key1} and \eqref{key2} are easily verified.
By Lemma \ref{belong}, $v$ is a solution of \eqref{floer} and hence
\eqref{key2} for $\widetilde{H_s}=H_r$.
By the definition of $\Sigma$, $v=\lim_{k\to \infty} \tau_k \cdot v_k$
and so 
$$\AC^j_{H_r}(v) =  \overleftarrow{\AC}^{-s_j^-}_{s_j^-}(v)= \lim_{k \to \infty} \overleftarrow{\AC}^{-s_j^-}_{s_j^-}(\tau_k \cdot v_k).$$
To prove that $v$ satisfies \eqref{key1}, it then suffices to show that
\begin{equation*}
\overleftarrow{E}^{-s_j^-}_{s_j^-}(\tau_k \cdot v_k)-\|H_k\|^- \leq \overleftarrow{\AC}^{-s_j^-}_{s_j^-}(\tau_k \cdot v_k)\leq \|H_k\|^+ -\epsilon
\end{equation*}
 for sufficiently large $k$.
The proof of these inequalities is entirely similar to the proofs
of \eqref{jk} and \eqref{k estimate} and is left to the reader.

We still must show that $v$ is a right asymptotic 
spanning disc for some nonconstant periodic orbit $y$ in $\PP(H_r)$.
Since $v$ satisfies 
\eqref{floer} and has energy less than $r(M,\om)$, there
is a sequence $s_j^- \to -\infty$ such that 
 $\lim_{j \to \infty} v(s^-_j,t) = y(t)\in \PP(H_r).$ 
We now prove that there is a sequence $s_j^+ \to +\infty$ such that 
$v(s^+_j,t)$ converges to a constant map. 

\begin{Lemma}\label{C minus}
A fixed point of the $\R$-action which belongs to $\VV_c \subset \CC$ is a
constant periodic orbit of $H_r$ contained in $C$.
\end{Lemma}

\begin{proof}
A fixed point of the $\R$-action in $\VV_c$ is an element of 
$\PP(H_r)$ which gets mapped by the evaluation to $U_{\scriptscriptstyle{(2r/3-\delta, +\infty)}}$.
Any periodic orbit of $H_r$ which enters $U_{\scriptscriptstyle{(2r/3-\delta, +\infty)}}$ must be 
constant. 
\end{proof}

Let $\tau_j \to +\infty$ be a sequence of 
positive numbers.
Passing to a subsequence, if necessary, we may assume that $\tau_j\cdot v$
converges to a map $\hat{v}$ in $\Sigma$. Since $v$ is not in $C$, it is 
nonconstant by Lemma \ref{C minus}. Lemma \ref{nonincrease} then implies 
that $\hat{v}$ belongs to $\VV_c$.

The limit $\hat{v}$ is also a fixed point of the $\R$-action. To prove this, we consider
the function $\tau \mapsto \overleftarrow{\AC}^0_{-\infty}(\tau \cdot v)$. 
Lemma \ref{action bounds} implies that this function 
is bounded from below by $-B$. Since $v$ is nonconstant, Lemma \ref{nonincrease} implies that it is also strictly decreasing. Hence, the limit
$$
\lim_{\tau \to + \infty}  \overleftarrow{\AC}^0_{-\infty}(\tau \cdot v)=d
$$
for some $d \geq -B$. By continuity, we then have
$$
\overleftarrow{\AC}^0_{-\infty}(\tau \cdot \hat{v})=\lim_{j \to + \infty}\overleftarrow{\AC}^0_{-\infty}((\tau+ \tau_j) \cdot v)=  \lim_{j \to + \infty}\overleftarrow{\AC}^0_{-\infty}(\tau_j \cdot v)=
\overleftarrow{\AC}^0_{-\infty}(\hat{v})
$$
for every $\tau$. Thus, $\hat{v}$ is an element of $\PP(H_r)$ by Lemma \ref{nonincrease}.

It now follows from Lemma \ref{C minus} that $\hat{v}$ is a constant periodic 
orbit of $H_r$ corresponding to some point $p$ in 
$U_{\scriptscriptstyle{(2r/3-\delta, +\infty)}}.$ Since the sequence $\tau_j \cdot v$ 
converges to the constant map $p$ in the $C^{\infty}_{loc}$-topology on $\CC$, the maps 
 $\tau_j \cdot v(0,t) = v(\tau_j,t)$ converge to $p$ in $C^{\infty}(S^1,M)$. 
Setting $s^+_j = \tau_j$ we have verified that $v$ is a right
asymptotic spanning disc for $y$.

Finally, as in \S \ref{case 1}, the fact that $y$ is nonconstant 
follows easily from \eqref{key1} and the fact that $v$ is nonconstant.


\begin{thebibliography}{FHS}


\bibitem[Al]{al}
P. Albers, A note on local Floer homology, Preprint 2006;
math.SG/0606600.





\bibitem[Bo1]{bo1}
 P. Bolle, Une condition de contact pour les sous-vari\'et\'es
 co\"isotropes d'une vari\'et\'e symplectique,
 \emph{C. R. Acad.\ Sci.\ Paris, S\'erie I}, \textbf{322} (1996), 83--86.

 \bibitem[Bo2]{bo2}
 P. Bolle, A contact condition for p-dimensional submanifolds of a
 symplectic manifold ($2\leq p\leq n$), \emph{Math.\ Z.}, \textbf{227}
 (1998), 211--230.

 \bibitem[Ch]{ch}
 Y. Chekanov, Lagrangian intersections, symplectic energy, and areas of
 holomorphic curves,  \emph{Duke Math.\ J.},
 \textbf{95} (1998), 213--226.

 \bibitem[Dr]{dr}
 D. Dragnev, Symplectic rigidity, symplectic fixed points and global perturbations of
 Hamiltonian systems, Preprint 2005, math.SG/0512109.



 \bibitem[EH]{eh}
 I. Ekeland, H. Hofer, Two symplectic fixed-point theorems with applications to Hamiltonian
 dynamics, \emph{J. Math. Pures Appl.} \textbf{68} (1989), no.4, 467–489 (1990).







  


  \bibitem[FHS]{fhs}
  A. Floer, H. Hofer, D. Salamon,
  Transversality in elliptic Morse theory for the symplectic action,
  \emph{Duke Math.\ J.}, \textbf{80} (1995), 251--292.

 
\bibitem[Gi]{gi:coisotropic}
 V.L. Ginzburg, Coisotropic intersections, to appear in \emph{Duke Math. J.}





\bibitem[Ho]{ho1}
H. Hofer, On the topological properties of symplectic maps,
\emph{Proc.\ Royal Soc.\ Edinburgh}, \textbf{115} (1990), 25--38.




 \bibitem[HZ]{hz:book} H. Hofer, E. Zehnder,
 \emph{Symplectic Invariants and Hamiltonian Dynamics}, Birk\"auser, 1994.






 \bibitem[Ke]{ke2}
 E. Kerman, Hofer's geometry and Floer theory under the quantum limit, Preprint 2007; math.SG/0703064.










%


\bibitem[Mo]{mo}
J. Moser, A fixed point theorem in symplectic geometry, \emph{Acta Math.}, \textbf{141} (1978), 17–34.









 \bibitem[PSS]{pss}
 S. Piunikhin, D. Salamon, M. Schwarz,
 Symplectic Floer--Donaldson theory and quantum cohomology,
 in \emph{Contact and symplectic geometry (Cambridge,
 1994)}, 171--201, Publ.\ Newton Inst., 8,
 Cambridge Univ.\ Press, Cambridge, 1996.


\bibitem[Po]{po}
L. Polterovich, Symplectic displacement energy for Lagrangian submanifolds, \emph{Ergod. Th. $\&$ Dynam. Sys.}, \textbf{13}, (1993), pp.\ 357--367.




 \bibitem[Sa]{Sa}
 D.A. Salamon, Lectures on Floer homology, in
 \emph{Symplectic Geometry and Topology},  Eds: Y. Eliashberg and
 L. Traynor, IAS/Park City Mathematics series, \textbf{7}, 1999,
 pp.\ 143--230.


  \bibitem[Sc]{schl}
  F. Schlenk,
  Applications of Hofer's geometry to Hamiltonian dynamics, \emph{Comment.\ Math.\ Helv.}, \textbf{81} (2006) 105--121.













 \bibitem[Vi]{vi:torus}
 C. Viterbo, A new obstruction to embedding Lagrangian tori, \emph{Invent. Math.}, \textbf{100} (1990), 301--320.


 

\end{thebibliography}
\end{document}